\documentclass{anstrans}
\title{Coarse Mesh Iteration Approach for Analytical 1D Multigroup $S_N$ Eigenvalue Problems}
\author{Jilang Miao,$^{*1}$ and Miaomiao Jin$^{*}$}

\institute{
$^{*}$Department of Nuclear Engineering, The Pennsylvania State University, University Park, 16802 PA, USA
}
\footnotetext[1]{\footnotesize Corresponding author: Jilang Miao (jlmiao@psu.edu)}

\usepackage{graphicx} 
\usepackage{subfigure}

\usepackage{booktabs} 
\usepackage{microtype} 
\usepackage{upgreek}
\usepackage{amsmath}
\usepackage{bbm}
\usepackage{algorithm}
\usepackage{algpseudocode}

\newcommand{\SN}{S$_N$~}
\DeclareMathOperator*{\argmin}{arg\,min}

\expandafter\def\expandafter\normalsize\expandafter{%
    \normalsize%
    \setlength\abovedisplayskip{1pt}%
    \setlength\belowdisplayskip{1pt}%
    \setlength\abovedisplayshortskip{1pt}%
    \setlength\belowdisplayshortskip{1pt}%
}

\begin{document}
\section{Introduction}
The discrete ordinates method, commonly referred to as the \SN method~\cite{Car1953}, involves discretizing the particle transport equation in its differential form. The determination of the particle fluxes is based on the straightforward evaluation of the transport equation at a limited number of discrete angular directions, or ordinates. Furthermore, quadrature relationships are utilized to replace integrals over angles, simplifying calculations through summations over these discrete ordinates~\cite{hebert2009applied}. 

Previous research introduced an accurate eigenvalue solver for multigroup \SN equations in slab geometry, employing the analytical \SN solution for each homogeneous subregion~\cite{analytical2A2G,AMGSN2023}. The method characterizes the solution in each subregion through an expansion based on the eigensystem determined by neutron cross sections in the material. The expansion coefficients are obtained by solving a linear system that incorporates continuity conditions at the interfaces and boundary conditions of the angular fluxes. The eigenvalue is determined by seeking the root of the determinant of the boundary condition matrix. Additionally, an analytical fixed source solver was developed and applied in the power iteration to address eigenvalue problems~\cite{AMGSNfixedsource}. Despite achieving a significant speedup compared to the sweeping-based \SN method, the full potential of the analytical \SN method was not harnessed when the source term was represented by piece-wise constant functions on a fine mesh.

In this study, we formulate an iterative approach for addressing the eigenvalue problem on coarse meshes, building upon the fixed source capability introduced in~\cite{AMGSNfixedsource}. The source term is expressed through an eigensystem expansion on the identical coarse mesh as the eigenfunctions. The power iteration process continually refines the eigensystem and expansion coefficients until convergence is achieved. We subsequently illustrate that the coarse mesh iteration method attains comparable high accuracy to the analytical fixed source solver on fine meshes while realizing a substantial speedup.
\section{Theory}
\subsection{Formulation of Fixed Source Problem}

For a given number of energy groups, denoted as $g=1,...,G$, 
and a quadrature set $\left . \{\mu_n,\omega_n\} \right | _{n=1,...,N}$, 
the transport equation for the angular flux $\psi_{g,n}$ is expressed in Eq~\ref{eq::sn1d}. 

\begin{equation}
\begin{aligned}
& \mu_n \frac{\partial}{\partial x} \psi_{g,n}(x)+  \Sigma_{t, g} \psi_{g,n}(x)= 
   \sum_{n^{\prime}, g^{\prime}} \omega_{n \prime} \Sigma_{s, g^{\prime} n^{\prime} \rightarrow g n} \psi_{g^{\prime}, n^{\prime}}(x) \\
&+\frac{1}{k_{eff}}  \sum_{n^{\prime}, g^{\prime}} \omega_{n^{\prime}}  \nu \Sigma_{f, g^{\prime} n^{\prime} \rightarrow g n} \psi_{g^{\prime}, n^{\prime}}(x) 
\label{eq::sn1d}
\end{aligned}
\end{equation}

The angular flux $\psi_{g,n}$ can be compactly aggregated in a vector $\Psi(x)$ of length $NG$. 
This vector consists of $G$ blocks, each having a length of $N$. For a specific block $g$ ($g=1,...,G$), 
it contains the angular fluxes $\left . \psi_{g,n} \right|_{n=1,\cdots,N}$. 
Consequently, we can denote $\psi_{g,n}(x)$ as $\Psi_{g N + n}(x)$.

With the same convention as in \cite{AMGSN2023} to organize the cross-sections and quadrature sets into matrices,
Eq~\ref{eq::sn1d} can be written in matrix form as,
\begin{equation}
    \partial_x \Psi(x) = A \Psi(x) 
    \label{eq::sn1d matrix}
\end{equation}
where
\begin{equation}
  A = \frac{F}{k_{eff}} + S - T
  \label{eq::TFStoA}
\end{equation}
Herein, $F$, $S$, and $T$ are the respective fission, scattering, and total cross sections multiplied by \SN quadrature set parameters  $\left . \{\mu_n,\omega_n\} \right | _{n=1,...,N}$, as in Eq~\ref{eq::sn1d}.

For eigenvalue problems, we can split matrix $A$ into two terms ($A=A_0+F_n$). Hence, Eq. \ref{eq::sn1d matrix} becomes, 
\begin{equation}
    \partial_x \Psi(x) = A_0 \Psi(x) + F_n(k_{eff}) \Psi(x) 
    \label{eq::sn1d A0Fn}
\end{equation}
Then, $F_n(k_{eff}) \Psi(x) $ part is viewed as external source $Q$. 
\begin{equation}
    \partial_x \Psi(x) = A_0 \Psi(x) + Q(k_{eff},x)
    \label{eq::sn1d A0FnQ}
\end{equation}
With power iteration at step $n$, we solve a fixed source problem as in Eq.~\ref{eq::sn1d A0Fn iter},
\begin{equation}
    \partial_x \Psi^{(n+1)}(x) = A_0 \Psi^{(n+1)}(x) + Q^{(n)}(k_{eff}^{(n)},x)
    \label{eq::sn1d A0Fn iter}
\end{equation}
In general cases, we have 
\begin{align}
     A_0 &= -T+S \\
     F_n &= \frac{F}{k_{eff}}
    \label{eq::def A0Fn}
\end{align}
To apply Wielandt's shift of $k_e$ ~\cite{brown2007wielandt}, we can have
\begin{align}
     A_0 &= -T+S + \frac{F}{k_{e}}\\
     F_n &= \frac{F}{k_{eff}} - \frac{F}{k_{e}}
\end{align}

\subsection{Coarse Mesh Iteration}
Within this section, we derive the equations essential for updating the source, eigenfunction, and eigenvalue during the iteration process on a coarse mesh. Both the eigenvalue and the source term are expressed as vectors of eigensystem expansion coefficients on the coarse mesh. The primary challenges involve evaluating the integral for the fixed source problem and updating the coefficients following changes in the eigensystem due to the updated $k_{eff}$.

\subsubsection{Representation of the source term}
With the eigenvalue problem represented by Eq~\ref{eq::sn1d matrix},
the angular flux can be found by solving Eq.~\ref{eq::Psi=PG beta} for coefficients in $\beta$. 
\begin{equation}
    \Psi(x) = P \Gamma(x) \beta
    \label{eq::Psi=PG beta}
\end{equation}
where matrices $P$ and $\Gamma(x)$ can be constructed from the eigensystem of matrix $A$~\cite{AMGSN2023}. 

As the iteration proceeds, $A$ is updated due to its dependency on eigenvalue $k_{eff}^{(n)}$,
so will be $P$ and $\Gamma(x)$.
Therefore, 
\begin{equation}
    \Psi^{(n)}(x) = P^{(n)} \Gamma^{(n)}(x) \beta^{(n)} + c^{(n)}
    \label{eq::Psi=PG beta n}
\end{equation}
The constant vector $c^{(n)}$ is added for flexibility, 
which should converge to $\mathbf{0}$. 
With angular flux given in Eq.~\ref{eq::Psi=PG beta n},
the source term in Eq~\ref{eq::sn1d A0Fn iter} can be explicitly written as 

\begin{equation}
Q^{(n)}(x) = F_n^{(n)} \left(  P^{(n)} \Gamma^{(n)}(x) \beta^{(n)} + c^{(n)} \right)
\label{eq::Qn iter}
\end{equation}
Here, we emphasize the dependency of $F_n$ on $k_{eff}^{(n)}$ by the superscript $(n)$.

\subsubsection{Solution of the angular flux}
Given Eq.~\ref{eq::Qn iter}, the angular flux for next iteration $(n+1)$ can be solved as Eq.~\ref{eq::psi n1}~\cite{AMGSNfixedsource}.
\begin{equation}
    \Psi^{(n+1)}(x) = P_0\Gamma_0 (x) \left (  \alpha^{(n+1)} + 
    \int_{x_0}^x d\xi \Gamma_0 (-\xi) P_0 ^{-1} Q^{(n)}(\xi)
    \right)
    \label{eq::psi n1}
\end{equation}
where $P_0$ and $\Gamma_0(x)$ are constructed from eigensystem of matrix $A_0$ in Eq~\ref{eq::sn1d A0Fn iter},
and $\alpha^{(n+1)}$ are to-be-determined coefficients~\cite{AMGSN2023}. 
Because $\Gamma_0(x)$ and $\Gamma^{(n)}(x)$ are either exponential or trigonometric functions,
the integral in Eq.~\ref{eq::psi n1} can be found analytically.
If the eigenvalues of  $A_0$ are all real, $\Gamma_0(x)$ is diagonal.

For ease of notation, 
we define the following matrices, 
\begin{align}
     F^{\dagger (n)} &= P_0^{-1} F^{(n)}_n P^{(n)} \\
     F_c^{(n)} &= P_0^{-1} F^{(n)}_n  \\
 B^{(n)}&=  {P^{(n)}}^{-1} A^{(n)} P^{(n)} \\
 B_0 &=  {P_0 }^{-1} A_0 P_0 
                 \label{eq::FdnFcn}
\end{align}
After derivation (skipped in this summary due to length limit), the  integral in  Eq.~\ref{eq::psi n1}  can be evaluated as
\begin{equation}
\begin{aligned}
  &  \int_{x_0}^x d\xi \Gamma_0 (-\xi) P_0 ^{-1} Q^{(n)}(\xi) \\
  &= 
\Gamma_0(-x) F^{\ddagger (n)} \Gamma^{(n)}(x) 
-\Gamma_0(-x_0) F^{\ddagger (n)} \Gamma^{(n)}(x_0)  \\
& -B_0^{-1}\left( \Gamma_0(-x) - \Gamma_0(-x_0) \right) F^{(n)}_c c^{(n)}
    \label{eq::psi n1 inte} 
\end{aligned}
\end{equation}
where
\begin{equation}
F^{\ddagger (n)} = U * \left(  F^{\dagger (n)} {B^{(n)}}^T -B_0 F^{\dagger (n)}\right)
\label{eq::Fd2Fdd}
\end{equation}
Here, $*$ is element-wise matrix product,
and
\begin{equation}
U_{i,j} = \frac{1}{\left| \lambda_{0,i} - \lambda^{(n)}_j \right|^2}
\label{eq::defU}
\end{equation}
where $\lambda_{0,i}$ and  $ \lambda^{(n)}_j$ are the $i$th and $j$th eigenvalues of $A_0$ and $A^{(n)}$, respectively. 
With the analytical representation of the source integral in Eq.~\ref{eq::psi n1},
$\alpha^{(n+1)}$ can now be solved from boundary conditions from a linear system as discussed in ~\cite{AMGSN2023,AMGSNfixedsource}. 

\subsubsection{Evaluation of the eigenvalue $k_{eff}$}
We define the fission source at iteration $(n)$ as the sum of new neutrons from fission across all energy groups, all discrete angles and over the spatial space,
\begin{equation}
\left|\left| F\Psi^{(n)} \right|\right| = \sum_{g,n} \omega_{n} \int dx \left( F \Psi^{(n)}(x) \right) _{g,n}
\label{eq::defsumQ}
\end{equation}
The angular flux is normalized by
\begin{equation}
  \left|\left| F\Psi^{(n)} \right|\right|  = 1
  \label{eq::iter norm}
\end{equation}
Therefore,  the eigenvalue can be updated as 
\begin{equation}
k^{(n+1)}_{eff} =   \left|\left| F\Psi^{(n+1)} \right|\right|  
\label{eq::kn1}
\end{equation}

The integral to obtain $k^{(n+1)}_{eff}$ is reduced to the integral of the angular flux vector $\Psi^{(n+1)}(x)$
 (with $F$ being constant within each region in the coarse mesh) and can be calculated as
\begin{equation}
\begin{aligned}
  &  \int_{x_L}^{x_R} dx  \Psi^{(n+1)}(x) =  \\
  & 
    P_0 B_0^{-1} \left(  \Gamma_0(x_R) - \Gamma_0(x_L)\right) \\
  & \left[ \alpha^{(n+1)} - \Gamma_0(-x_L)F^{\ddagger (n)}\Gamma^{(n)}(x_L) \beta^{(n)}  + B_0^{-1} \Gamma_0(-x_L) F_c^{(n)} c^{(n)} \right] \\
 & + P_0 F^{\ddagger (n)} {B^{(n)}}^{-1} \left( \Gamma^{(n)}(x_R) - \Gamma^{(n)}(x_L)\right)  \beta^{(n)}  \\
 & - P_0 B_0^{-1} F_c^{(n)} c^{(n)} (x_R-x_L)
   \label{eq::intePsin1}
\end{aligned}
\end{equation}
where $x_L$ and $x_R$ are the positions of left and right boundaries of the region.

\subsubsection{Projection of angular flux on the new eigensystem}
With $k_{eff}$ updated in Eq~\ref{eq::kn1},
and the resultant new matrices $A^{(n+1)}$, $P^{(n+1)}$, $B^{(n+1)}$, $\Gamma^{(n+1)}(x)$, 
we then need to find the new coefficients $\beta^{(n+1)}$ and $c^{(n+1)}$ in Eq~\ref{eq::Psi=PG beta n}.
Eq~\ref{eq::Psi=PG beta n} can be viewed as the projection of $\Psi$ on the new eigensystem ($P\Gamma$). To solve the projection coefficients ($\beta$ and $c$), Eq~\ref{eq::Psi=PG beta n} can be formulated as an optimization problem.
We first evaluate $\Psi (x)$ on $J$ spatial points ($\left. \{\Psi^{(n+1)}(x_j)\} \right |_{j=1,\cdots,J}$) according to Eq.~\ref{eq::psi n1}.
Then the expansion coefficients ($\beta$ and $c$) can be identified by minimizing the square error between the angular flux before and after the projection step, i.e., 
\begin{equation}
  \begin{aligned}
    &\beta^{(n+1)},c^{(n+1)}  = \\
    &\argmin_{\{\beta,c \}} \sum_{j=1}^J \sum_{g,n} \omega_n  
    \left| P^{(n+1)} \Gamma^{(n+1)}(x_j) \beta + c   - \Psi^{(n+1)}_{g,n}(x_j)  \right|^2\\
    &= \argmin_{\{\beta,c \}} \sum_{j=1}^J \sum_{i}^{NG} \Omega_i  
    \left| M_2(x_j) \beta + c   - \Psi^{(n+1)}_{i}(x_j)  \right|^2
      \label{eq::projection opt}
\end{aligned}
    \end{equation}
where vector $\Omega$ and matrix $M_2$ are defined as the following to simplify the notation.
\begin{equation}
\Omega_{gN+n} = \omega_n 
  \label{eq::defWw}
\end{equation}

\begin{equation}
M_2(x_j) = P^{(n+1)}\Gamma^{(n+1)}(x_j) 
\label{eq::defM2}
\end{equation}

Similar to the derivation of linear regression~\cite{freedman2009statistical}, we can find the coefficients as the solution of Eq.~\ref{eq::projection sln}. 
\begin{equation}
  \begin{split}
    \left[
    \begin{array}{cc}
      W & W_0 \\
      \bar{M_2} & \mathbbm{1}
    \end{array}
      \right]
\left[\begin{array}{c}
\beta^{(n+1)} \\
c^{(n+1)}
\end{array} \right]
= 
\left[\begin{array}{c}
V \\
\bar{\psi}
\end{array} \right]
  \end{split}
\label{eq::projection sln}
\end{equation}
where
\begin{equation}
  W_{m,k}  = \sum_{i,j} \Omega_j {M_2(x_i)}_{j,k} {M_2(x_i)}_{j,m} 
  \label{eq::def W}
\end{equation}

\begin{equation}
  {W_0}_{m,j}  = \sum_{i} \Omega_j  {M_2(x_i)}_{j,m}   
\label{eq::def W0}
\end{equation}

\begin{equation}
 {\bar{M_2}}_{m,k}  = \sum_{i}   {M_2(x_i)}_{m,k}   
\label{eq::def M2b}
\end{equation}

\begin{equation}
  {V}_{m}  = \sum_{i,j} \Omega_j \Psi_j(x_i)  {M_2(x_i)}_{j,m}   
\label{eq::def V}
\end{equation}

\begin{equation}
 \bar{\psi}_{m}  = \sum_{i}   \Psi_m(x_i)
\label{eq::def psibar}
\end{equation}
Finally, the source term of iteration $(n+1)$ is updated according to Eq.~\ref{eq::Qn iter}. 
\subsubsection{Algorithm for coarse mesh iteration}
Based on the developed theory,
the algorithm to perform the eigenvalue power iteration is summarized in Alg.~\ref{alg::coarse iter}. 
\begin{algorithm}
  \caption{Eigenvalue power iteration with analytical multigroup \SN solver}
  \label{alg::coarse iter}
\begin{algorithmic}
\State initialize $k^{(0)}_{eff}$ 
\For{each distinct material}
\State find block-diagonalization matrices $P_0$ and $B_0$ for $A_0$
\State find block-diagonalization matrices $P^{(0)}$ and $B^{(0)}$ for $A^{(0)}=A_0 + F_n^{(0)}$
\EndFor
\State initialize angular flux $\left. \{\Psi^{(0)}(x_j)\} \right |_{j=1,\cdots,J}$ for $J$ points
\State find $\beta^{(0)},c^{(0)}$  (Eq.~\ref{eq::projection sln})
\While{ error metric above threshold }
    \State solve coefficients $\alpha^{(n+1)}$ (Eq.~\ref{eq::psi n1})
    \State update $k^{(n+1)}_{eff}$ from $\alpha^{(n+1)}$, $\beta^{(n)}$ and $c^{(n)}$ (Eq~\ref{eq::kn1})
\For{each distinct material}
\State find block-diagonalization matrices $P^{(n+1)}$ and $B^{(n+1)}$ for $A^{(n+1)}=A_0 + F_n^{(n+1)}$
\EndFor
\State evaluate angular flux $\left. \{\Psi^{(n+1)}(x_j)\} \right |_{j=1,\cdots,J}$ for $J$ points
\State find $\beta^{(n+1)},c^{(n+1)}$  (Eq.~\ref{eq::projection sln})
\State calculate error metric such as norm of $\Psi^{(n+1)}-\Psi^{(n)}$
\EndWhile
\end{algorithmic}
\end{algorithm}

Notably, for each region, $J$ can be rather sparse to solve for the eigensystem expansion coefficients. Also, the angular flux evaluated at these $J$ spatial points can be used to calculate the error metric as noted in Alg.~\ref{alg::coarse iter}. 

Note that the equations here are derived for a homogeneous problem. However, it can be extended to heterogeneous problems, with additional consideration of the interface continuity conditions as discussed in ~\cite{AMGSN2023,AMGSNfixedsource}. Our case study in the next section will demonstrate this capability.

\section{Results}
As a test case, we study a 35 cm slab with 3 regions.
The reactor core is located within [-15 cm, 15 cm].
The reflector is within [-17.5 cm, -15 cm] and [15 cm, 17.5cm].
The system has vacuum boundary condition on both ends.
Two-group cross-sections  for the core and reflector materials are generated with OpenMC~\cite{romano2013openmc,boyd2019multigroup} for a typical fuel pincell, which can be found in~\cite{AMGSNfixedsource}. 

Alg.~\ref{alg::coarse iter} is run for \SN in the cases of $N=2,~4,~8,~16$ with Gauss-Legendre quadrature sets.
The system is divided into a coarse mesh of size $R=24$ ($2$ in each reflector and $20$ in the core).
With each grid, a uniform mesh of size $J=2$ is chosen, where the angular fluxes at the $J$ mesh centers are evaluated to calculate the eigensystem expansion coefficients (Eq.~\ref{eq::projection sln}).
The initial guess of the angular flux is assumed to be isotropic and vary according to $\Psi_{g,n}(x) \propto |x|$ .
The iteration is terminated when the $L^2$ norm of angular flux change between two consecutive generations is below the threshold in Eq.~\ref{eq::exerr}.
\begin{equation}
||\Psi^{(n)}-\Psi^{(n-1)}||_2 < 10 ^{-6} \sqrt{\frac{R J N}{700}}
\label{eq::exerr}
\end{equation}

We will compare the efficiency of the coarse mesh iteration approach developed in this work
and the fine mesh iteration approach developed in~\cite{AMGSNfixedsource}.
Note that the scaling factor $\sqrt{\frac{R J N}{700}}$ in Eq.~\ref{eq::exerr} is intentionally used so that this criterion is equivalent to termination criterion ($10^{-6}$) for scalar flux on a mesh of size $700$ in ~\cite{AMGSNfixedsource}.

\subsection{Accuracy of the coarse mesh iteration method}
A reference solution is generated using OpenMC~\cite{romano2013openmc} multigroup mode with
the same geometric configuration, boundary conditions and cross-sections as the test case. 
The simulation tracks $10^6$ neutrons per generation.
The neutrons are simulated for $200$ inactive generations and tallies are collected for the next $800$ active generations to compute  
scalar fluxes, angular fluxes and $k_{eff}$.
The fluxes are tallied on a spatially uniform mesh of size $700$ for each energy group.
In addition, the angular fluxes are tallied over a specific polar angle range corresponding to the \SN quadrature set.

The normalization step (Eq.\ref{eq::iter norm}) makes the fluxes from \SN directly comparable with MC. 
Table~\ref{tab::keff} shows the $k_{eff}$ from OpenMC and the different orders of the analytical \SN solvers.
It clearly shows how higher-order solution approaches the MC reference.

\begingroup
\scriptsize
\begin{table}[htb]
  \centering
  \caption{Computed $k_{eff}$ on a coarse mesh compared with MC reference.}
  \begin{tabular}{llr}\toprule
Method      & $k_{eff}$      & $k_{eff}$ - $k_{eff,MC}$ (pcm) 
\\ \midrule
MC reference & 1.24953 $\pm$ 0.00002 &  \\
\midrule
Analytical $S_2$        & 1.24737 & -216 \\
Analytical $S_4$        & 1.24936 & -17 \\
Analytical $S_8$        & 1.24949 & -4 \\
Analytical $S_{16}$     & 1.24952 & -1 \\
    \bottomrule
\end{tabular}
  \label{tab::keff}
\end{table}
\endgroup

In addition, Table ~\ref{tab::keff coarse fine} indicates the $k_{eff}$ of all orders from the coarse mesh match those
from the fine mesh of size $700$ in~\cite{AMGSNfixedsource} with less than  1 pcm difference.

\begingroup
\scriptsize
\begin{table}[htb]
  \centering
  \caption{$k_{eff}$ from analytical methods on coarse mesh vs fine mesh.}
  \begin{tabular}{lllr}\toprule
\SN order      & coarse mesh     & fine mesh  & $k_{eff}$ difference (pcm)  \\
$S_2$        & 1.247372 & 1.247371 & 0.1 \\
$S_4$        & 1.249365 & 1.249364 & 0.1 \\
$S_8$        & 1.249490 & 1.249490 & 0.0 \\
$S_{16}$     & 1.249519 & 1.249518 & 0.1 \\
    \bottomrule
\end{tabular}
  \label{tab::keff coarse fine}
\end{table}
\endgroup

With the converged solution represented by coefficients on the coarse mesh of size $24$,
we can then evaluate the angular flux at any location.
They are evaluated at the same mesh of size $700$ for pointwise comparison with MC. We observe a significant enhancement in accuracy with increasing \SN orders. In the case of scalar fluxes, the point-wise relative error diminishes from approximately $10\%$ in $S_2$ to about $0.1\%$ in $S_{16}$. Similarly, for angular fluxes ($\omega_n \psi_{g,n}(x)$), the point-wise relative error decreases from around $30\%$ in $S_2$ to approximately $0.5\%$ in $S_{16}$. The plots of scalar fluxes, angular fluxes, and point-wise relative error with MC reference closely resemble those in \cite{AMGSNfixedsource} and are therefore omitted in this summary.

\subsection{Efficiency of the coarse mesh iteration method}
In this demonstration, we highlight the efficiency advantages of the analytical coarse mesh iteration method.
Fig. 1(a)
plots the $L^2$ norm of flux changes versus the number of iterations.
It shows that both the coarse mesh and fine mesh methods converge at the same rate at all the \SN orders.
They all converge with the equivalent criteria after around $25$ iterations.
This suggests that different orders of \SN methods exhibit similar dominance ratios, despite differences in $k_{eff}$. 
Fig. 1(a)
also shows that with the Wielandt's shift $k_e=1.4$, 
the solver is significantly accelerated and converges within $10$ iterations. 
Note that the assumption of real eigenvalues of $A_0$ to simplify the integral in Eq.~\ref{eq::psi n1} is justified
for the problem without Wielandt's shift or $k_e$ above 1.4. 
For lower $k_e$ such as $1.35$, $A_0$ would have complex eigenvalues, necessitating modifications to Eq.~\ref{eq::psi n1 inte}. This will be addressed in future work.

We then assess the computational cost associated with each iteration. 
Fig. 1(b) and Fig. 1(c)
plot the $L^2$ norm of flux change versus time,
measured in the unit of the average time required to solve one iteration in the case of the coarse mesh analytical $S_{16}$ without Wielandt's shift. 
Fig. 1(b)
shows that the coarse mesh method is significantly faster than the fine mesh method.
With the same $k_{eff}$ solution,
the coarse mesh method has $18\times$, $9\times$, $14\times$, $3\times$ speed up for the $S_2$, $S_4$, $S_8$, $S_{16}$ orders, respectively.
The efficiency disparity between the coarse mesh and fine mesh iteration methods stems from two competing factors.
In the fine mesh method, despite fluxes being represented by coefficients on a coarse mesh, they must be evaluated at fine mesh locations to update the source term, which assumes a piece-wise constant nature and necessitates the fine mesh. 
The increased cost associated with evaluating flux values makes the fine mesh method less efficient.
By contrast, the coarse mesh method represents both flux and source with coefficients on the coarse mesh, but additional costs are incurred in updating the eigensystem for each material after each iteration.

Fig. 1(c)
further substantiates the speedup attributed to Wielandt's shift. 
It is noted that the average time per iteration remains the same after enabling Wielandt's shift,
This outcome aligns with expectations, as Wielandt's shift does not modify Alg.~\ref{alg::coarse iter} and the tested $k_e$ values do not introduce complexity from complex eigenvalues. 

\begin{figure}[htbp]
\centering
\includegraphics[width=0.43\textwidth]{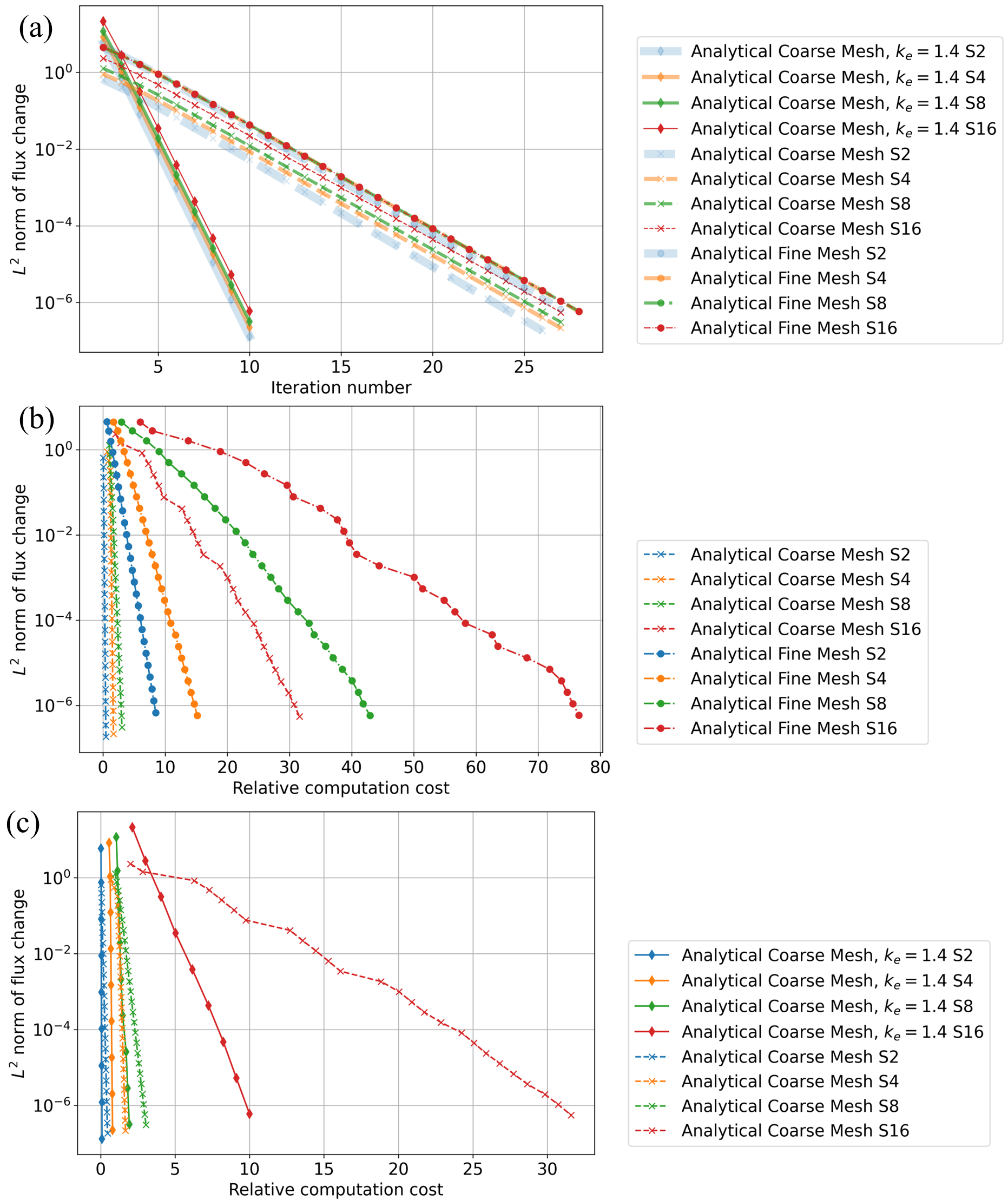}
\caption{Convergence of fine and coarse mesh analytical \SN methods. (a) Norm of flux change vs. iteration number. (b) Computation time compared between coarse mesh and fine mesh method. (c) Computation time compared between coarse mesh method with and without Wielandt's shift. }
\label{fig::conv}
\end{figure}

\section{Conclusions}
In this study, we introduced an efficient iteration method based on a coarse mesh to solve the eigenvalue problem of analytical multigroup \SN equations in slab geometry. For the slab problem homogenized from a typical pincell, we achieved an eigenvalue accuracy of $216$ pcm for the $S_2$ solution and $1$ pcm for the $S_{16}$ solution. We also observed high accuracy in scalar and angular fluxes. In comparison to the fine mesh iteration method, which necessitates a piecewise constant source term, the analytical coarse mesh iteration method demonstrates a substantial speedup while maintaining the same high level of accuracy.

\section{Acknowledgments}
This work is supported by the Department of Nuclear Engineering, The Pennsylvania State University.


\footnotesize
\bibliographystyle{ans}
\bibliography{bibliography}

\begin{thebibliography}{1}
\newcommand{\enquote}[1]{``#1''}

\bibitem{Car1953}
\MakeUppercase{B.~G. Carlson}, \enquote{Solution of the Transport Equation by
  {S$_n$} Approximations,} Tech. Rep. LA-1599, Los Alamos Scientific Laboratory
  (1953).

\bibitem{hebert2009applied}
\MakeUppercase{A.~H{\'e}bert}, \emph{Applied Reactor Physics}, Presses inter
  Polytechnique (2009).

\bibitem{analytical2A2G}
\MakeUppercase{J.~Miao} and \MakeUppercase{M.~Jin}, \enquote{An Analytic Method
  for Solving Static Two-group, {1D} Neutron Transport Equations,}
  \emph{Transactions of the American Nuclear Society}, \textbf{127}, 1068--1071
  (2022).

\bibitem{AMGSN2023}
\MakeUppercase{J.~Miao} and \MakeUppercase{M.~Jin}, \enquote{An Accurate S$_N$
  Method for Solving Static Multigroup Neutron Transport Equations in Slab
  Geometry,} \emph{Transactions of the American Nuclear Society}, \textbf{129},
  926--929 (2023).

\bibitem{AMGSNfixedsource}
\MakeUppercase{J.~Miao} and \MakeUppercase{M.~Jin}, \enquote{Developing an
  Analytical Fixed Source Solver for the 1D Multigroup $S_N$ Equations,}
  \emph{arXiv identifier to be updated} (2024).

\bibitem{brown2007wielandt}
\MakeUppercase{F.~Brown} \MakeUppercase{et~al.}, \enquote{Wielandt acceleration
  for MCNP5 Monte Carlo eigenvalue calculations,} in \enquote{Joint
  International Topical Meeting on Mathematics \& Computation and
  Supercomputing in Nuclear Applications (M\&C+ SNA 2007), Monterey,
  California,}  (2007).

\bibitem{freedman2009statistical}
\MakeUppercase{D.~A. Freedman}, \emph{Statistical models: theory and practice},
  cambridge university press (2009).

\bibitem{romano2013openmc}
\MakeUppercase{P.~K. Romano} and \MakeUppercase{B.~Forget}, \enquote{The
  {OpenMC} monte carlo particle transport code,} \emph{Annals of Nuclear
  Energy}, \textbf{51}, 274--281 (2013).

\bibitem{boyd2019multigroup}
\MakeUppercase{W.~Boyd}, \MakeUppercase{A.~Nelson}, \MakeUppercase{P.~K.
  Romano}, \MakeUppercase{S.~Shaner}, \MakeUppercase{B.~Forget}, and
  \MakeUppercase{K.~Smith}, \enquote{Multigroup cross-section generation with
  the OpenMC Monte Carlo particle transport code,} \emph{Nuclear Technology},
  \textbf{205}, \emph{7}, 928--944 (2019).

\end{thebibliography}
\end{document}